\documentclass[12pt, a4paper]{article}

\usepackage{epsfig}
\usepackage{amssymb}
\usepackage{amsthm}%% The amsthm package provides extended theorem environments
\usepackage{amsmath}
\usepackage{float}
\usepackage{multicol}
\usepackage{multirow}
\usepackage{cite}
\usepackage{enumerate}
\usepackage{xcolor}
\usepackage{tabu}
\def\ms{\medskip}
\def\nt{\noindent}

\definecolor{vividviolet}{rgb}{0.62, 0.0, 1.0}

\def\Z{\mathbb Z}
\newtheoremstyle{de}%name
  {10pt}          % space above
  {10pt}  % space below
  {\rm}  % bofy font
  {}%{\parindent}     % ident - empty=no indent,  \parindent= paragraph indent
  {\bf}  % thm head font
  {. }    % punctuation after thm head
  { }    % space after thm head: `` ``=normal \newline=linebreak
  {}     % thm head specification
\theoremstyle{de}

\newtheorem{example}{Example}[section]

\newtheoremstyle{theorem}%name
  {10pt}          % space above
  {10pt}  % space below
  {\it}  % bofy font
  {}%{\parindent}     % ident - empty=no indent,  \parindent= paragraph indent
  {\bf}  % thm head font
  {. }    % punctuation after thm head
  { }    % space after thm head: `` ``=normal \newline=linebreak
  {}     % thm head specification
\theoremstyle{theorem}

\newtheorem{theorem}{Theorem}[section]
\newtheorem{corollary}[theorem]{Corollary}%[section]

\numberwithin{equation}{section}
\def\Z{\mathbb{Z}}
\def\N{\mathbb{N}}

\usepackage[mathlines]{lineno}
\modulolinenumbers[2]
\newcommand*\patchAmsMathEnvironmentForLineno[1]{%
\expandafter\let\csname old#1\expandafter\endcsname\csname #1\endcsname  \expandafter\let\csname oldend#1\expandafter\endcsname\csname end#1\endcsname  \renewenvironment{#1}%
{\linenomath\csname old#1\endcsname}%
{\csname oldend#1\endcsname\endlinenomath}}%
\newcommand*\patchBothAmsMathEnvironmentsForLineno[1]{%
\patchAmsMathEnvironmentForLineno{#1}%
\patchAmsMathEnvironmentForLineno{#1*}}%
\AtBeginDocument{%
\patchBothAmsMathEnvironmentsForLineno{equation}%
\patchBothAmsMathEnvironmentsForLineno{align}%
\patchBothAmsMathEnvironmentsForLineno{flalign}%
\patchBothAmsMathEnvironmentsForLineno{alignat}%
\patchBothAmsMathEnvironmentsForLineno{gather}%
\patchBothAmsMathEnvironmentsForLineno{multline}%
}

\begin{document}
%\baselineskip18truept
%\normalsize
\begin{center}
{\mathversion{bold}\Large \bf New Families of tripartite graphs with local antimagic chromatic number 3}

\bigskip
{\large  Gee-Choon Lau$^a$, Wai Chee Shiu$^b$ }\\

\medskip

\emph{{$^a$}77D, Jalan Suboh, 85000 Segamat, Johor, Malaysia}\\
%\emph{Universiti Teknologi MARA (Segamat Campus),}\\
%\emph{85000, Johor, Malaysia.}\\
\emph{geeclau@yahoo.com}\\

\medskip
\emph{{$^b$}Department of Mathematics,}\\
\emph{The Chinese University of Hong Kong,}\\
\emph{Shatin, Hong Kong, P.R. China.}\\
\emph{wcshiu@associate.hkbu.edu.hk}\\

\end{center}

\begin{abstract}
For a graph $G(V,E)$ of size $q$, a bijection $f : E(G) \to [1,q]$ is a local antimagc labeling if it induces a vertex labeling $f^+ : V(G) \to \N$ such that $f^+(u) \ne f^+(v)$, where $f^+(u)$ is the sum of all the incident edge label(s) of $u$, for every edge $uv \in E(G)$. In this paper, we make use of matrices of fixed sizes to construct several families of infinitely many tripartite graphs with local antimagic chromatic number 3. % Consequently we are able to determine the local antimagic chromatic number of many infinite families of related graphs. In [Local antimagic vertex coloring of corona product graphs $P_n \circ P_k$, {\it Advances in Computer Science Research}, {\bf 96}, 10.2991/acsr.k.220202.014], the authors claimed that $\chi_{la}(P_n \circ P_3) = 6$ (Theorem 2.2) by showing that $\chi_{la}(P_n\circ P_3)\le 6$ in Lemma 2.2.1 and $\chi_{la}(P_n \circ P_3)\ne 5$ in Lemma 2.2.2. We then give counterexamples to odd $n\ge 3$ by giving a local antimagic 5-coloring of $P_n\circ P_3$.
\ms

\noindent Keywords: Local antimagic  chromatic number, tripartite, regular, disconnected

\noindent 2020 AMS Subject Classifications: 05C78; 05C69.
\end{abstract}

\baselineskip18truept
\normalsize

\section{Introduction}
Let $G=(V, E)$ be a connected graph of order $p$ and size $q$.
A bijection $f: E\rightarrow \{1, 2, \dots, q\}$ is called a \textit{local antimagic labeling}
if $f^{+}(u)\neq f^{+}(v)$ whenever $uv\in E$,
where $f^{+}(u)=\sum_{e\in E(u)}f(e)$ and $E(u)$ is the set of edges incident to $u$.
The mapping $f^{+}$ which is also denoted by $f^+_G$ is called a \textit{vertex labeling of $G$ induced by $f$}, and the labels assigned to vertices are called \textit{induced colors} under $f$.
The \textit{color number} of a local antimagic labeling $f$ is the number of distinct induced colors under $f$, denoted by $c(f)$.  Moreover, $f$ is called a \textit{local antimagic $c(f)$-coloring} and $G$ is {\it local antimagic $c(f)$-colorable}. The \textit{local antimagic chromatic number} $\chi_{la}(G)$ is defined to be the minimum number of colors taken over all colorings of $G$ induced by local antimagic labelings of $G$~\cite{Arumugam}. Let $G+H$ and $mG$ denote the disjoint union of graphs $G$ and $H$, and $m$ copies of $G$, respectively. For integers $c < d$, let $[c,d] = \{n\in\Z\;|\; c\le n\le d\}$. Very few results on the local antimagic chromatic number of regular graphs are known (see~\cite{Arumugam, LauLiShiu}).
Throughout this paper, we let $V(aP_2\vee O_m) = \{u_i, v_i, x_j\;|\; 1\le i\le a, 1\le j\le m\}$ and $E(aP_2 \vee O_m) = \{u_ix_j, v_ix_j, u_iv_i\;|\; 1\le i\le a, 1\le j\le m\}$. We also let $V(a(P_2\vee O_m)) = \{u_i, v_i, x_{i,j}\;|\; 1\le i\le a, 1\le j\le m\}$ and $E(a(P_2\vee O_m)) = \{u_ix_{i,j}, v_ix_{i,j}, u_iv_i\;|\; 1\le i\le a, 1\le j\le m\}$.

\ms \nt In~\cite{Haslegrave}, the author proved that all connected graphs without a $P_2$ component admit a local antimagic labeling. Thus, $O_m, m\ge 1$ and $aP_2, a\ge 1$ are the only families of regular graphs without local antimagic chromatic number.  In~\cite{Arumugam}, it was shown that $\chi_{la}(aP_2\vee O_1) = 3$ for $a\ge 1$. In the following sections, we extend the ideas in~\cite{LauShiu-join, LSPN-odd} to construct various families of tripartite graphs of size $(4n+1)\times (2k+1)$ and $(4n+3)\times (2k+1)$, for $n, k\ge 1$, respectively, and proceed to prove that all these graphs have local antimagic chromatic number 3.

%\newpage
\section{Graphs of size $(4n+1)\times (2k+1)$}

\nt For $k\ge 1$, we now consider the following $(4n+1)\times (2k+1)$ matrix for $2\le j\le n$.  Note that when $n=1$, the required $5\times (2k+1)$ matrix is given by rows $f(u_i, x_{i,1})$, $f(u_i, x_{i,2})$, $f(u_iv_i)$, $f(v_ix_{i,1})$ and $f(v_1x_{i,2})$ of the matrix below. Moreover, the entries in column $k+1$ appears in both parts of the matrix.

\[\fontsize{7}{10}\selectfont
\begin{tabu}{|c|[1pt]c|c|c|c|c|c|[1pt]c|[1pt]c}\hline
i & 1 & 2 & 3 & \cdots & k-1 & k & k+1 &  \\\tabucline[1pt]{-}
\multirow{2}{1.2cm}{$f(u_ix_{i,1})$} & k+2+  & k+3 +  & k+4+ & \multirow{2}{0.3cm}{$\cdots$} & 2k+   & 2k+1+ & 1+   \\
 & n(8k+4)& n(8k+4) & n(8k+4) & & n(8k+4) & n(8k+4) & n(8k+4) \\ \hline
\multirow{2}{1.2cm}{$f(u_ix_{i,2})$} & -2k-2+ & -2k-4 + & -2k-6+ & \multirow{2}{0.3cm}{$\cdots$} &  -4k+2+ &  -4k & - 2k-1+    \\
 & n(8k+4)& n(8k+4) & n(8k+4) &  & n(8k+4) & n(8k+4) & n(8k+4) \\ \hline
\vdots & \vdots & \vdots & \vdots & \cdots & \vdots & \vdots & \vdots \\\hline
\multirow{2}{1.5cm}{$f(u_ix_{i,2j-1})$} & 9k+6 & 9k+7 & 9k+8 & \multirow{2}{0.3cm}{$\cdots$}  &  10k+4  & 10k+5 & 8k+5    \\
 & (n-j)(8k+4)& (n-j)(8k+4) & (n-j)(8k+4) &  & (n-j)(8k+4) & (n-j)(8k+4) & (n-j)(8k+4) \\ \hline
\multirow{2}{1.2cm}{$f(u_ix_{i,2j})$} &  5k+2 & 5k+1 & 5k & \multirow{2}{0.3cm}{$\cdots$} & 4k+4 &     4k+3 & 6k+3 \\
 & (n-j)(8k+4)& (n-j)(8k+4) & (n-j)(8k+4) &  & (n-j)(8k+4) & (n-j)(8k+4) & (n-j)(8k+4) \\ \hline
\vdots & \vdots & \vdots & \vdots & \cdots & \vdots & \vdots & \vdots & \\\hline
%f(u_ix_{i,2n-1}) & 9k+6+ & 9k+7+ & 9k+8+ & \cdots & 10k+4+ & 10k+5+ & 8k+5+  \\\hline
%f(u_ix_{i,2n}) & 5k+2+ & 5k+1+ & 5k + & \cdots & 4k+4 & 4k+3 & 6k+3 +   \\\hline
(u_iv_i) & 1 & 2 & 3 & \cdots & k-1 & k & k+1 & \\\hline
f(v_ix_{i,1}) & 3k+2 & 3k+3 & 3k+4 & \cdots & 4k & 4k+1 & 4k+2  \\\hline
f(v_ix_{i,2}) &   8k+4 & 8k+2 & 8k & \cdots & 6k+8 &    6k+6 & 6k+4  \\\hline
\vdots & \vdots & \vdots & \vdots & \cdots & \vdots & \vdots & \vdots & \\\hline
\multirow{2}{1.5cm}{$f(v_ix_{i,2j-1})$} & -5k-2+ & -5k-1+  & -5k+ & \multirow{2}{0.3cm}{$\cdots$} & -4k-4+  &  -4k-3+ & -4k-2 +    \\
& j(8k+4) & j(8k+4) & j(8k+4) &  & j(8k+4) & j(8k+4) & j(8k+4)   \\\hline
\multirow{2}{1.2cm}{$f(v_ix_{i,2j})$} & -k+ & -k-1+  & -k-2+ & \multirow{2}{0.3cm}{$\cdots$} & -2k+3+ &    -2k+1+ & -2k+    \\
& j(8k+4) & j(8k+4) & j(8k+4) &  & j(8k+4) & j(8k+4) & j(8k+4)   \\\hline
\vdots & \vdots & \vdots & \vdots & \cdots & \vdots & \vdots & \vdots & \\\hline
\end{tabu}
\]
\[\fontsize{7}{10}\selectfont
\begin{tabu}{|c|[1pt]c|[1pt]c|c|c|c|c|c|}\hline
i  & k+1 & k+2 & k+3 & \cdots & 2k-1 & 2k & 2k+1 \\\tabucline[1pt]{-}
\multirow{2}{1.2cm}{$f(u_ix_{i,1})$} & 1+  & 2+ & 3+ & \multirow{2}{0.3cm}{$\cdots$} & k-1+ & k+ & k+1+   \\
 & n(8k+4)& n(8k+4) & n(8k+4) &   & n(8k+4) & n(8k+4) & n(8k+4) \\  \hline
\multirow{2}{1.2cm}{$f(u_ix_{i,2})$} & -2k-1+ & -2k-3+ & -2k-5+ & \multirow{2}{0.3cm}{$\cdots$} & -4k+3+ & -4k+1+ & -4k-1+    \\
& n(8k+4)& n(8k+4) & n(8k+4) &  & n(8k+4) & n(8k+4) & n(8k+4) \\ \hline
\vdots & \vdots & \vdots & \vdots & \cdots & \vdots & \vdots & \vdots \\\hline
\multirow{2}{1.5cm}{$f(u_ix_{i,2j-1})$} & 8k+5  & 8k+6 & 8k+7 & \multirow{2}{0.3cm}{$\cdots$} & 9k+3 & 9k+4 & 9k+5 \\
 & (n-j)(8k+4)& (n-j)(8k+4) & (n-j)(8k+4) &   & (n-j)(8k+4) & (n-j)(8k+4) & (n-j)(8k+4) \\ \hline
\multirow{2}{1.2cm}{$f(u_ix_{i,2j})$} &  6k+3 & 6k+2 & 6k+1 & \multirow{2}{0.3cm}{$\cdots$} & 5k+5 & 5k+4 & 5k+3 \\
 & (n-j)(8k+4)& (n-j)(8k+4) & (n-j)(8k+4) &   & (n-j)(8k+4) & (n-j)(8k+4) & (n-j)(8k+4) \\ \hline
\vdots & \vdots & \vdots & \vdots & \cdots & \vdots & \vdots & \vdots \\\hline
f(u_iv_i) & k+1 & k+2 & k+3 & \cdots & 2k-1 & 2k & 2k+1 \\\hline
f(v_ix_{i,1}) & 4k+2 &  2k+2 & 2k+3 & \cdots & 3k-1 & 3k & 3k+1 \\\hline
f(v_ix_{i,2}) & 6k+4   & 8k+3 & 8k+1 & \cdots & 6k+9 & 6k+7 & 6k+5\\\hline
 \vdots & \vdots & \vdots & \vdots & \cdots & \vdots & \vdots & \vdots \\\hline
\multirow{2}{1.5cm}{$f(v_ix_{i,2j-1})$} & -4k-2+  & -6k-2+ & -6k-1+ & \multirow{2}{0.3cm}{$\cdots$} & -5k-5+ & -5k-4+ & -5k-3+   \\
& j(8k+4)& j(8k+4) & j(8k+4) &  & j(8k+4) & j(8k+4) & j(8k+4) \\ \hline
\multirow{2}{1.2cm}{$f(v_ix_{i,2j})$} & -2k+   & 0+ & -1+ & \multirow{2}{0.3cm}{$\cdots$} & -k+3+ & -k+2+ & -k+1+    \\
 & j(8k+4)& j(8k+4) & j(8k+4) &  & j(8k+4) & j(8k+4) & j(8k+4) \\ \hline
\vdots & \vdots & \vdots & \vdots & \cdots & \vdots & \vdots & \vdots \\\hline
\end{tabu}
\]

\nt We now have the following observations.
\begin{enumerate}[(a)]
\item  For $n\ge 2$ and each $i\in [1,2k+1]$, the sum of the first $2n+1$ row entries is $f^+(u_i)= 2n(8k+4)-k+1+\sum\limits^{n}_{j=2} [2(n-j)(8k+4)+14k+8] = 8kn^2+6kn+4n^2+k+4n+1$.  Note that, this formula also holds when $n=1$.
\item For $n\ge 2$ and each $i\in [1,2k+1]$, the sum of the last $2n+1$ row entries is $f^+(v_i) = 11k+7+\sum\limits^{n}_{j=2} [2j(8k+4)-6k-2] = 8kn^2+2kn+4n^2+k+2n+1$.  Note that, this formula also holds when $n=1$.
\item For each $i\in [1,k]$ and $j\in [1,2n]$, each of $f(u_ix_{i,j}) + f(v_{2k+2-i}x_{2k+2-i,j})$, $f(v_ix_{i,j}) + f(u_{2k+2-i}x_{2k+2-i,j})$ and $f(u_{k+1}x_{k+1,j}) + f(v_{k+1}x_{k+1,j})$ is a constant $n(8k+4)+4k+3$.
\item Suppose $2k+1=(2r+1)(2s+1)$, $r,s\ge 1$.  For each $a \in [1,r]$ and $j\in [1,2n]$, each of
\begin{align}&\sum^{2s+1}_{b=1} [f(u_{(a-1)(2s+1)+b}x_{(a-1)(2s+1)+b,j}) + f(v_{2k+2-(a-1)(2s+1)-b}x_{2k+2-(a-1)(2s+1)-b,j})],\label{eq-d1}\\
&\sum^{2s+1}_{b=1} [f(v_{(a-1)(2s+1)+b}x_{(a-1)(2s+1)+b,j}) + f(u_{2k+2-(a-1)(2s+1)-b}x_{2k+2-(a-1)(2s+1)-b,j})],\label{eq-d2}\\
&\sum^{2s+1}_{b=1} [f(u_{r(2s+1)+b}x_{r(2s+1)+b,j}) + f(v_{2k+2-r(2s+1)-b}x_{2k+2-r(2s+1)-b,j})],\label{eq-d3}\end{align}
is a constant $(2s+1)[n(8k+4)+4k+3]$.
\end{enumerate}

\nt Consider $G = (2k+1)P_2 \vee O_{2n}$. By Observations (a) and (b) above, we can now define a bijection $f : E(G) \to [1, (4n+1)(2k+1)]$ according to the table above. Clearly, for $1\le i\le 2k+1$, $f^+(u_i) > f^+(v_i)$.

\ms\nt Now, for each $i\in [1,k]$ and $j\in [1,2n]$, first delete the edges $v_ix_{i,j}$ and $v_{2k+2-i}x_{2k+2-i,j}$, and then add the edges  $v_{2k+2-i}x_{i,j}$ and $v_ix_{2k+2-i,j}$ with labels $f(v_{2k+2-i}x_{2k+2-i,j})$ and $f(v_ix_{i,j})$, respectively. Finally, we rename $x_{i,j}$ by $y_{i,j}$ and $x_{2k+2-i,j}$ by $z_{i,j}$. We still denote this new labeling by $f$. By Observation (c),  $f^+(y_{i,j})=f^+(z_{i,j})=n(8k+4)+4k+3$. It is easy to verify that $f^+(u_i) \ne f^+(v_i) \ne f^+(y_{i,j})$ for all possible $n, k$. We denote the resulting graph by $G_{2n}(k+1)$. Note that  $G_{2n}(k+1)$ has $k+1$ components.

\begin{theorem}\label{thm-G2n(k+1)} For $n, k\ge 1$, we have $\chi_{la}(G_{2n}(k+1)) = 3$. \end{theorem}

\begin{proof} From the above discussion, we know that $G_{2n}(k+1)$ is a tripartite graph with $k+1$ components that admits a local antimagic 3-coloring. The theorem holds. \end{proof}

\begin{example} Consider $n = 2$ and $k=4$. We have the following table.
\[%\fontsize{7}{10}\selectfont
\begin{tabu}{|c|[1pt]c|c|c|c|[1pt]c|[1pt]c|c|c|c|}\hline
i & 1 & 2 & 3 & 4 & 5 & 6 & 7 & 8 & 9  \\\tabucline[1pt]{-}
f(u_ix_{i,1}) & 78 & 79 & 80 & 81 & 73 & 74 & 75 & 76 & 77  \\\hline
f(u_ix_{i,2}) & 62 & 60 & 58 & 56 & 63 & 61 & 59 & 57 & 55 \\\hline
f(u_ix_{i,3}) & 42 & 43 & 44 & 45 & 37 & 38 & 39 & 40 & 41 \\\hline
f(u_ix_{i,4})  & 22  & 21 &  20 & 19 & 27  & 26 & 25 & 24 & 23 \\\hline
f(u_iv_i) & 1 & 2 & 3 & 4  & 5 & 6 & 7 & 8 & 9  \\\hline
f(v_ix_{i,1}) & 14 & 15 &  16 & 17 &  18  & 10 & 11 & 12 & 13 \\\hline
f(v_ix_{i,2}) & 36 & 34 & 32 & 30 & 28  & 35 & 33 & 31  & 29  \\\hline
f(v_ix_{i,3})  & 50 & 51 & 52 & 53 & 54  & 46 & 47 & 48 & 49 \\\hline
f(v_ix_{i,4})  & 68  & 67 & 66 & 65 & 64 & 72 & 71 & 70 & 69 \\\hline
\end{tabu}\]
By the construction above Theorem~\ref{thm-G2n(k+1)}, we have the graph $G_4(5)$ as shown below.
\begin{figure}[H]
\centerline{\epsfig{file=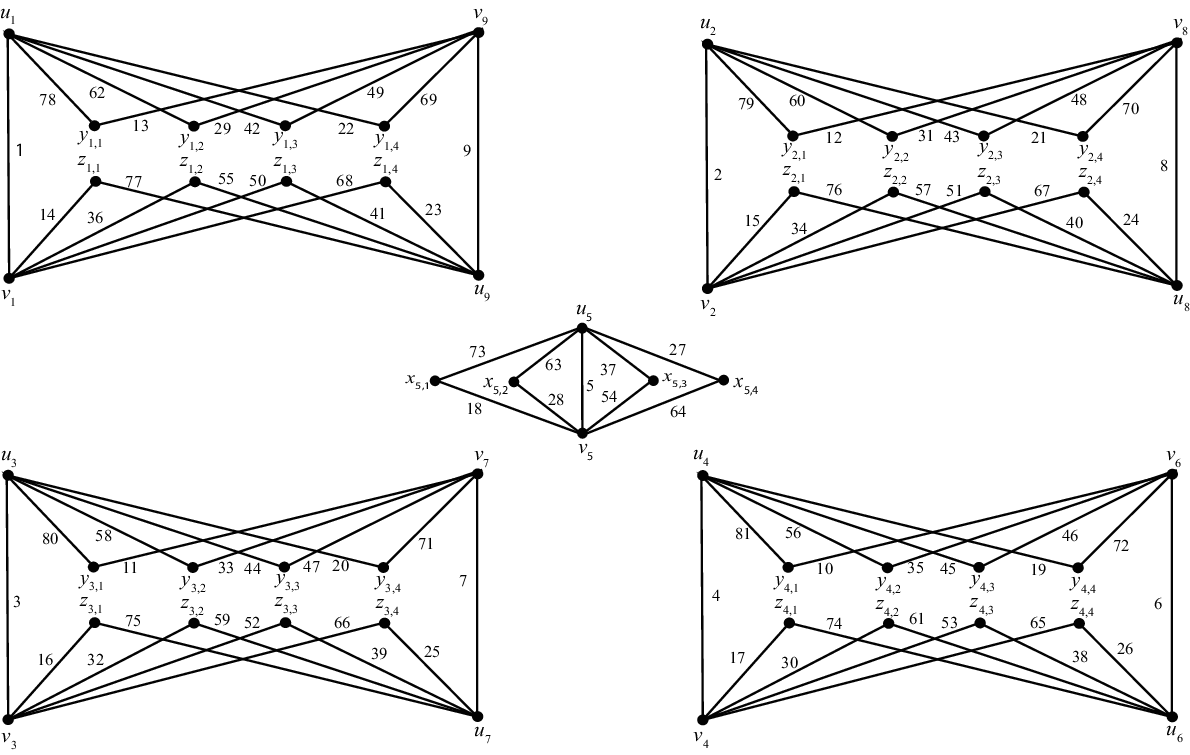, width=15cm}}
\caption{Graph $G_4(5)$.}\label{fig:G4(5))}
\end{figure}
\end{example}

\nt We may make use of Observation (d) to construct a new graph with local antimagic chromatic number 3 from $G_{2n}(k+1)$. Let us show an example first. Suppose $2k+1= (2r+1)(2s+1)$, $r, s\ge 1$.

\begin{example} Consider $n = 2$, $k=4$ again. Now we have $r=s=1$. Consider the graph $G=G_{2n}(k+1)$. Now $V(G)=\{u_i, v_i\;|\; 1\le i\le 9\}\cup\{y_{i,j}, z_{i,j}\;|\; 1\le i\le 4, 1\le j\le 4\}$.  From Observation (d) we have
\begin{align*}
f^+(y_{1,j})+f^+(y_{2,j})+f^+(y_{3,j}) &=[f(u_1x_{1,j})+f(v_9x_{9,j})]+[f(u_2x_{2,j})+f(v_8x_{8,j})]\\& \quad +[f(u_3x_{3,j})+f(v_7x_{7,j})]=273,\\
f^+(z_{1,j})+f^+(z_{2,j})+f^+(z_{3,j}) & =[f(v_1x_{1,j})+f(u_9x_{9,j})]+[f(v_2x_{2,j})+f(u_8x_{8,j})]\\& \quad+[f(v_3x_{3,j})+f(u_7x_{7,j})]=273,\\
f^+(y_{4,j})+f^+(x_{5,j})+f^+(z_{4,j}) &=[f(u_4x_{1,j})+f(v_6x_{2,j})]+[f(u_5x_{5,j})+f(v_5x_{5,j})]\\& \quad+[f(u_6x_{6,j})+f(v_4x_{4,j})]=273.
\end{align*}
For each $j\in[1,4]$, we (i) merge the vertices $y_{1,j}, y_{2,j}, y_{3,j}$ as a new vertex (still denote by $y_{1,j}$) of degree 6; (ii) merge the vertices $z_{1,j}, z_{2,j}, z_{3,j}$ as a new vertex (still denote by $z_{1,j}$) of degree 6; and (iii) merge $y_{4,j}, x_{5,j}, z_{4,j}$ (denote by $x_{5,j}$) of degree 6. %{\magenta for consistent, it is better to denote by $x_{5,j}$. If we do so, then we must revise the vertex names of the figure of $G_4(3,3)$.}

\begin{figure}[H]
\centerline{\epsfig{file=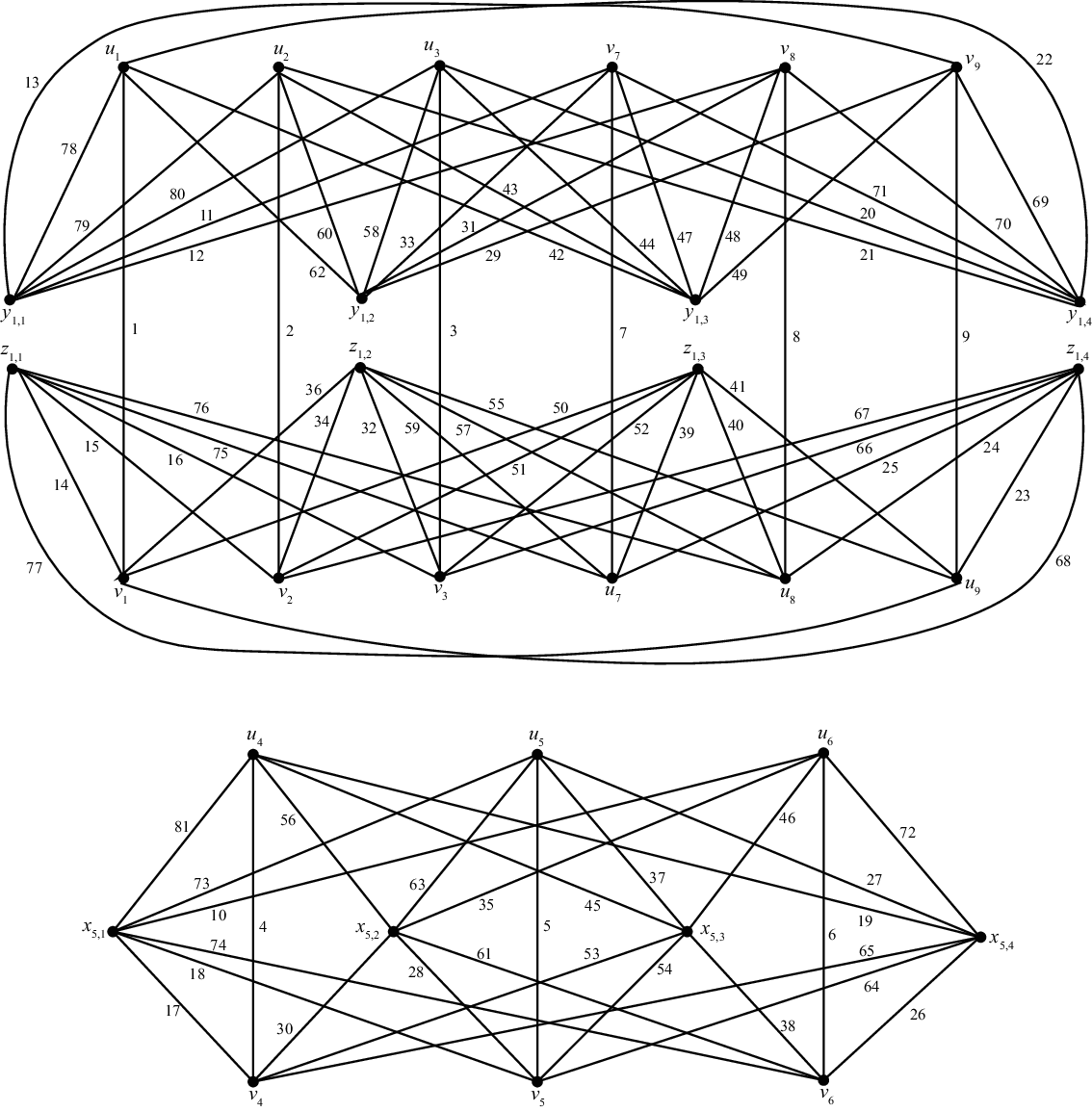, width=13cm}}
\caption{Graph $G_4(3,3)$.}\label{fig:G4(3,3))}  % {\magenta According to my notation introduced below, please change $x_{4,j}$ to $x_{5,j}$.}
\end{figure}

\end{example}

\ms\nt Suppose $2k+1 = (2r+1)(2s+1)$, $r,s\ge 1$. Consider the graph $G_{2n}(k+1)$. For each $a\in[1,r]$ and $j\in[1,2n]$, we can merge all $2s+1$ vertices in $\{y_{(a-1)(2s+1)+b, j}\;|\;b\in[1, 2s+1]\}$, $\{z_{(a-1)(2s+1)+b, j}\;|\;b\in[1, 2s+1]\}$, and $\{x_{r(2s+1)+b, j}\;|\;b\in[1, 2s+1]\}$. The new vertices are denoted by $y_{(a-1)(2s+1)+1,j}$, $z_{(a-1)(2s+1)+1,j}$ and $x_{k+1,j}$, respectively.
By equations~\eqref{eq-d1}, \eqref{eq-d2} and \eqref{eq-d3}, we have $f^+(y_{(a-1)(2s+1)+1,j})=f^+(z_{(a-1)(2s+1)+1,j})=f^+(x_{k+1,j})=(2s+1)[n(8k+4)+4k+3]$.  Let the graph just obtained be $G_{2n}(2r+1,2s+1)$.  Note that $G_{2n}(2r+1,2s+1)$ has $r+1$ components.

\begin{theorem}\label{thm-G2n(2r+1,2s+1)} For $n,r,s\ge 1$, we have $\chi_{la}(G_{2n}(2r+1,2s+1)) = 3$. \end{theorem}

\begin{proof} From the above discussion, we know that $2k+1 = (2r+1)(2s+1)$, $r,s\ge 1$ and $G_{2n}(2r+1,2s+1)$ is a tripartite graph with $r+1$ components that admits a bijective edge labeling $f$ with induced vertex labels $(1) =  (2s+1)[n(8k+4)+4k+3]$, $(2) = 8kn^2+6kn+4n^2+k+4n+1$, and $(3) = 8kn^2+2kn+4n^2+k+2n+1$. Clearly, $(2) > (3)$. We now show that $(1)\ne (2), (3)$. Now,
\begin{align*}
(1) - (2) & =  16kns-8kn^2+2kn+8ks-4n^2+8ns+3k+6s+2 \\
 &= (8kn+4n+3)(2s-n) + 2kn+8ks+3k+3n+2\\
 & >  0 \quad \mbox{ if } 2s\ge n.
\end{align*}
Otherwise, $2s -n \le - 1$ (equivalently, $-n\le -2s-1$), $(1) - (2) \le -6kn  - n -1 + 8ks + 3k = -n(6k+1) - 1 + 8ks + 3k \le (-2s-1)(6k+1) - 1 + 8ks + 3k = - 4ks-3k-2s-2 < 0$.  Thus, $(1) \ne (2)$. Similarly,
\begin{align*}
(1)-(3) & =  16kns-8kn^2+6kn+8ks-4n^2+8ns+3k+2n+6s+2\\
 &= (8kn+4n+3)(2s-n) + 6kn + 8ks + 3k + 5n + 2 \\
 & >  0  \quad \mbox{ if } 2s \ge n.
\end{align*}
If $2s-n = -1$, $(1)-(3) = -2kn-n-1+8ks+3k=-n(2k+1)-1+8ks+3k = (-2s-1)(2k+1)-1+8ks+3k = 4ks+k-2s-2 > 0$ since $k\ge 4$. Otherwise, $2s-n\le -2$ (equivalently, $-n\le -2s-2$), $(1) - (3) \le -10kn-3n-4+8ks+3k \le (-2s-2)(10k+3)-4+8ks+3k < 0$. Thus, $(1) \ne (3)$. Therefore, $f$ is a local antimagic 3-coloring. The theorem holds.
\end{proof}

%\nt Note that we can also apply the cut-merge process that gives us Theorem~\ref{thm-H2n(2r+1,2s+1)} to the graph Theorem~\ref{thm-G2n(2r+1,2s+1)} to get non-isomorphic (possibly connected) graphs also with local antimagic chromatic number 3.

\section{Graphs of size $(4n+3)\times (2k+1)$}

In what follows, we refer to the following $(4n+3)\times (2k+1)$ matrix to obtain results similar to Theorems~\ref{thm-G2n(k+1)}  and~\ref{thm-G2n(2r+1,2s+1)}. For $1\le j\le n$, we have

\[\fontsize{7}{10}\selectfont
\begin{tabu}{|c|[1pt]c|c|c|c|c|c|}\hline
i & 1 & 2 & 3 & \cdots & 2k & 2k+1  \\\tabucline[1pt]{-}
\vdots & \vdots & \vdots & \vdots & \cdots & \vdots & \vdots   \\\hline
\multirow{2}{1.5cm}{$f(u_ix_{i,2j-1})$} & 10k+5 + & 10k+4 +  & 10k+3  & \multirow{2}{0.3cm}{$\cdots$} & 8k+6 +& 8k+5 +  \\
 & (2n-j)(4k+2)  & (2n-j)(4k+2)  & (2n-j)(4k+2) &   & (2n-j)(4k+2)  & (2n-j)(4k+2)  \\\hline
\multirow{2}{1.2cm}{$f(u_ix_{i,2j})$} & 6k+4 + & 6k+5 + & 6k+6 + & \multirow{2}{0.3cm}{$\cdots$} & 8k+3 + & 8k+4 +    \\
 & (2n-j)(4k+2)  & (2n-j)(4k+2)  & (2n-j)(4k+2) &    & (2n-j)(4k+2)  & (2n-j)(4k+2)  \\\hline
\vdots & \vdots & \vdots & \vdots & \cdots & \vdots & \vdots   \\\hline
\multirow{2}{1.5cm}{$f(u_ix_{i,2n+1})$} &  2k+1 + & 2k+ & (2k-1)+  & \multirow{2}{0.3cm}{$\cdots$} & 2 + & 1 +  \\
 &  (n+1)(4k+2) & (n+1)(4k+2) & (n+1)(4k+2) &  & (n+1)(4k+2) & (n+1)(4k+2) \\\hline
f(u_iv_i) & 1 & 2 & 3 & \cdots & 2k & 2k+1  \\\hline
f(v_ix_{i,1}) & 4k+2 & 4k+1 & 4k & \cdots & 2k+3 & 2k+2   \\\hline
\vdots & \vdots & \vdots & \vdots & \cdots & \vdots & \vdots   \\\hline
\multirow{2}{1.2cm}{$f(v_ix_{i,2j})$} & 4k+3 + & 4k+4 +  & 4k+5 + &\multirow{2}{0.3cm}{$\cdots$} & 6k+2 + & 6k+3 +    \\
 & (j-1)(4k+2)  &  (j-1)(4k+2) &  (j-1)(4k+2) &    & (j-1)(4k+2) &  (j-1)(4k+2)  \\\hline
\multirow{2}{1.5cm}{$f(v_ix_{i,2j+1})$} & 8k + 4 + & 8k+3 +  & 8k+2 + & \multirow{2}{0.3cm}{$\cdots$} & 6k+5 + & 6k+4 +    \\
 & (j-1)(4k+2)  &  (j-1)(4k+2) &  (j-1)(4k+2) &   & (j-1)(4k+2)  & (j-1)(4k+2)    \\\hline
 \vdots & \vdots & \vdots & \vdots & \cdots & \vdots & \vdots   \\\hline
\end{tabu}
\]

\nt We now have the following observations. % {\Lau need to change (3) and (4)}
\begin{enumerate}[(1)]
\item For each column, the sum of the first $2n+2$ entries is $f^+(u_i) = (n+1)(3n+1)(4k+2) + n + 2k+2$.
\item For each column, the sum of the last $2n+2$ entries is $f^+(v_i) = (n+1)^2(4k+2)+n+1$.
\item For each $i\in[1,k]$ and $j\in[1,2n+1]$, each of $f(u_ix_{i,j}) + f(v_{2k+2-i}x_{2k+2-i,j})$, $f(v_ix_{i,j}) + f(u_{2k+2-i}x_{2k+2-i,j})$, and, $f(u_{k+1}x_{k+1,j}) + f(v_{k+1}x_{k+1,j})$ is a constant $(2n+2)(4k+2)+1$.
\item Suppose $2k+1=(2r+1)(2s+1)$, $r,s\ge 1$. For each $a\in[1,r]$ and $j\in[1,2n+1]$, each of
\begin{align}&\sum^{2s+1}_{b=1} [f(u_{(a-1)(2s+1)+b}x_{(a-1)(2s+1)+b,j}) + f(v_{2k+2-(a-1)(2s+1)-b}x_{2k+2-(a-1)(2s+1)-b,j})],\label{eq-d4}\\
&\sum^{2s+1}_{b=1} [f(v_{(a-1)(2s+1)+b}x_{(a-1)(2s+1)+b,j}) + f(u_{2k+2-(a-1)(2s+1)-b}x_{2k+2-(a-1)(2s+1)-b,j})],\label{eq-d5}\\
&\sum^{2s+1}_{b=1} [f(u_{r(2s+1)+b}x_{r(2s+1)+b,j}) + f(v_{2k+2-r(2s+1)-b}x_{2k+2-r(2s+1)-b,j})],\label{eq-d6}\end{align}
is a constant $(2s+1)[ (2n+2)(4k+2)+1]$.

\end{enumerate}

\nt Similar to graph $G_{2n}(k+1)$ in Theorem~\ref{thm-G2n(k+1)}, we also define $G_{2n+1}(k+1)$ of $k+1$ components similarly  such that the $i$-th component has vertex set $\{u_i, v_i, u_{2k+2-i}, v_{2k+2-i}, y_{i,j}, z_{i,j} \mid 1\le j\le 2n+1\}$ and edge set $\{u_iv_i, u_{2k+2-i}v_{2k+2-i}, u_iy_{i,j}, v_{2k+2-i}y_{i,j}, v_iz_{i,j}, u_{2k+2-i}z_{i,j}\mid 1\le j\le 2n+1\}$ for $1\le i\le k$, and the $(k+1)$-st component is the $P_2\vee O_{2n+1}$ with vertex set $\{u_{k+1}, v_{k+1}, x_{k+1,j} \mid 1\le j\le 2n+1\}$ and edge set $\{u_{k+1}v_{k+1}, u_{k+1}x_{k+1,j}, v_{k+1}x_{k+1,j}\mid 1\le j\le 2n+1\}$. Moreover, by Observation (3), $f^+(y_{i,j}) = f^+(z_{i,j}) = (2n+2)(4k+2)+1$. It is easy to verify that $f^+(u_i) \ne f^+(v_i) \ne f^+(y_{i,j})$ for all possible $n,k$.

\begin{theorem}\label{thm-G2n+1(k+1)} For $n, k\ge 1$, $\chi_{la}(G_{2n+1}(k+1)) = 3$.  \end{theorem}

\begin{proof} From the discussion above, we know $G_{2n+1}(k+1)$ is a tripartite graph with $k+1$ components that admits a local antimagic 3-coloring. The theorem holds.
\end{proof}

\nt For $2k+1 = (2r+1)(2s+1), r,s\ge 1$, by Observation (4) above, we also define $G_{2n+1}(2r+1,2s+1)$ as in Theorem~\ref{thm-G2n(2r+1,2s+1)} with $r+1$ components and similar vertex set with vertices $y_{(a-1)(2s+1)+1,j}$, $z_{(a-1)(2s+1)+1,j}$ and $x_{k+1,j}$ for $1\le a\le 2r+1$, $1\le j\le 2n+1$. By equations~\eqref{eq-d4}, \eqref{eq-d5} and \eqref{eq-d6}, we have  $f^+(y_{(a-1)(2s+1)+1,j})=f^+(z_{(a-1)(2s+1)+1,j})=f^+(x_{k+1,j})=(2s+1)[(2n+2)(4k+2)+1]$.

\begin{theorem}\label{thm-G2n+1(2r+1,2s+1)}   For $n, r ,s\ge 1$, we have $\chi_{la}(G_{2n+1}(2r+1,2s+1)) = 3$. \end{theorem}

\begin{proof} Similar to the proof of Theorem~\ref{thm-G2n(2r+1,2s+1)}, we know $2k+1 = (2r+1)(2s+1)$, $r,s\ge 1$ and $G_{2n+1}(2r+1,2s+1)$ is a tripartite graph with $r+1$ components that admits a bijective edge labeling $f$ with induced vertex labels $(1) = (2s+1)[(2n+2)(4k+2)+1]$, $(2) = (n+1)(2n+1)(4k+2)+n+2k+2$ and $(3) = (n+1)^2(4k+2)+n+1$. Clearly, $(2) > (3)$. We now show that $(1) \ne (2), (3)$.

\ms\nt Now,
\begin{align*}
(1) - (2) & = -8kn^2 +16kns-4kn+16ks-4n^2 +8ns+2k-3n+10s+1\\
   & = (8kn+4n+4k+5)(2s-n) + 2n + 8ks + 2k + 1\\
  & > 0 \quad \mbox{ if }  2s\ge n.
\end{align*}
If $2s-n \le -1$, $(1) -(2) \le -8kn -2n-2k-4+8ks \le (-2s-1)(8k+2) - 2k-4+8ks < 0$. Thus, $(1) \ne (2)$. Similarly,

\begin{align*}
(1) -(3) & = -4kn^2 +16kns+16ks-2n^2 +8ns+4k-n+10s+2\\
  & = (4kn+2n+2)(4s-n) + n + 16ks + 2s + 4k +  2 \\
 & > 0  \quad \mbox{ if } 4s\ge n.
\end{align*}
If $4s - n \le -1$, $(1)- (3) \le -4kn - n + 16ks + 2s + 4k \le (-4s-1)(4k+1) + 16ks+2s+4k = -2s -1 < 0$. Thus, $(1) \ne (3)$. Therefore, $f$ is a local antimagic 3-coloring. The theorem holds.
\end{proof}

\begin{example} Take $n=2$, $k=4$, we have the following table and graph $G_5(5)$ with the defined labeling.

\[
\begin{tabu}{|c|[1pt]c|c|c|c|c|c|c|c|c|}\hline
i & 1 & 2 & 3 & 4 & 5 & 6 & 7 & 8 & 9   \\\tabucline[1pt]{-}
f(u_ix_{i,1}) & 99 & 98 & 97 & 96 & 95 & 94 & 93 & 92 & 91  \\\hline
f(u_ix_{i,2}) & 82 & 83 & 84 & 85 & 86 & 87 & 88 & 89 & 90  \\\hline
f(u_ix_{i,3}) & 81 & 80 & 79 & 78 & 77 & 76 & 75 & 74 & 73   \\\hline
f(u_ix_{i,4}) & 64 & 65 & 66 & 67 & 68 & 69 & 70 & 71 & 72  \\\hline
f(u_ix_{i,5}) & 63 & 62 & 61 & 60 & 59 & 58 & 57 & 56 & 55 \\\hline
f(u_iv_i) & 1 & 2 & 3 & 4  & 5 & 6 & 7 & 8 & 9 \\\hline
f(v_ix_{i,1}) & 18 & 17 & 16 & 15 &  14 & 13 &  12 & 11 &  10   \\\hline
f(v_ix_{i,2}) & 19  & 20 & 21 & 22  & 23 & 24 & 25 & 26 & 27   \\\hline
f(v_ix_{i,3})  & 36 & 35 & 34 & 33  & 32 & 31 & 30 & 29 & 28  \\\hline
f(v_ix_{i,4})  & 37  & 38 & 39 & 40 & 41 & 42 & 43 & 44 & 45  \\\hline
f(v_ix_{i,5}) & 54 & 53 & 52 & 51  & 50 & 49 & 48 & 47 & 46 \\\hline
\end{tabu}\]

\begin{figure}[H]
\centerline{\epsfig{file=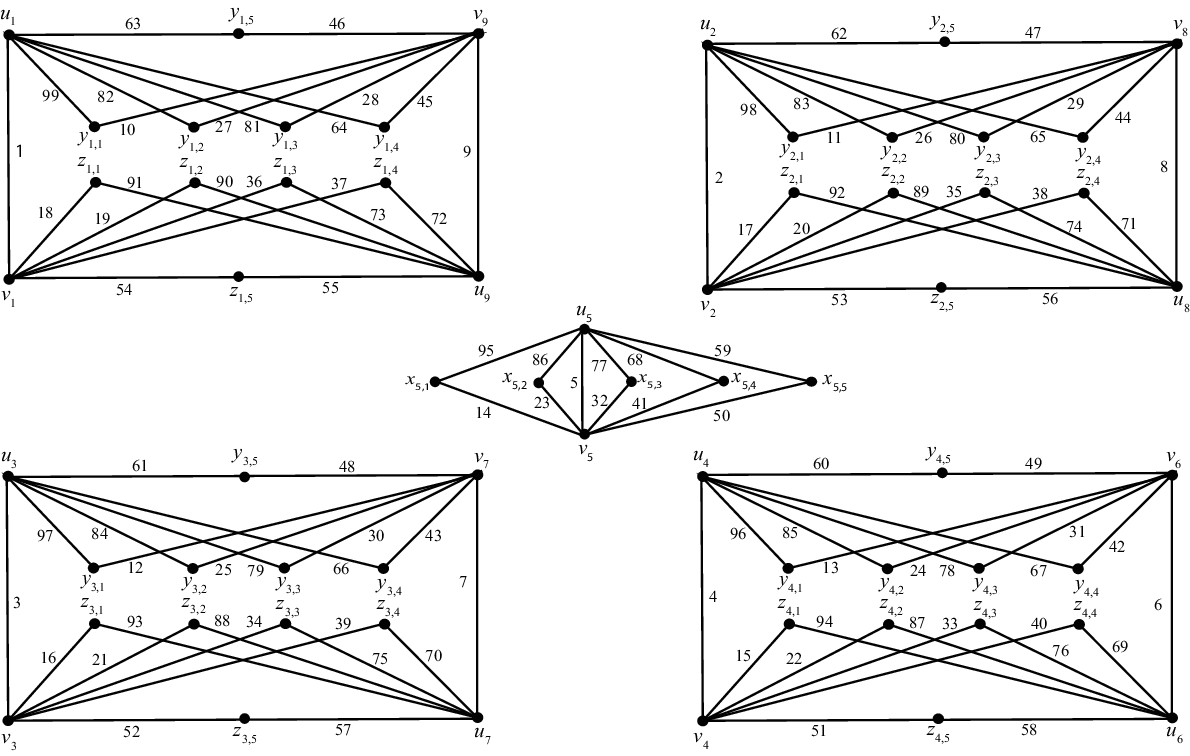, width=14cm}}
\caption{Graph $G_5(5)$.}\label{fig:G5(5))}
\end{figure}

\ms\nt If we take $r=s=1$, we can get $G_5(3,3)$ which is a 6-regular graph.
\begin{figure}[H]
\centerline{\epsfig{file=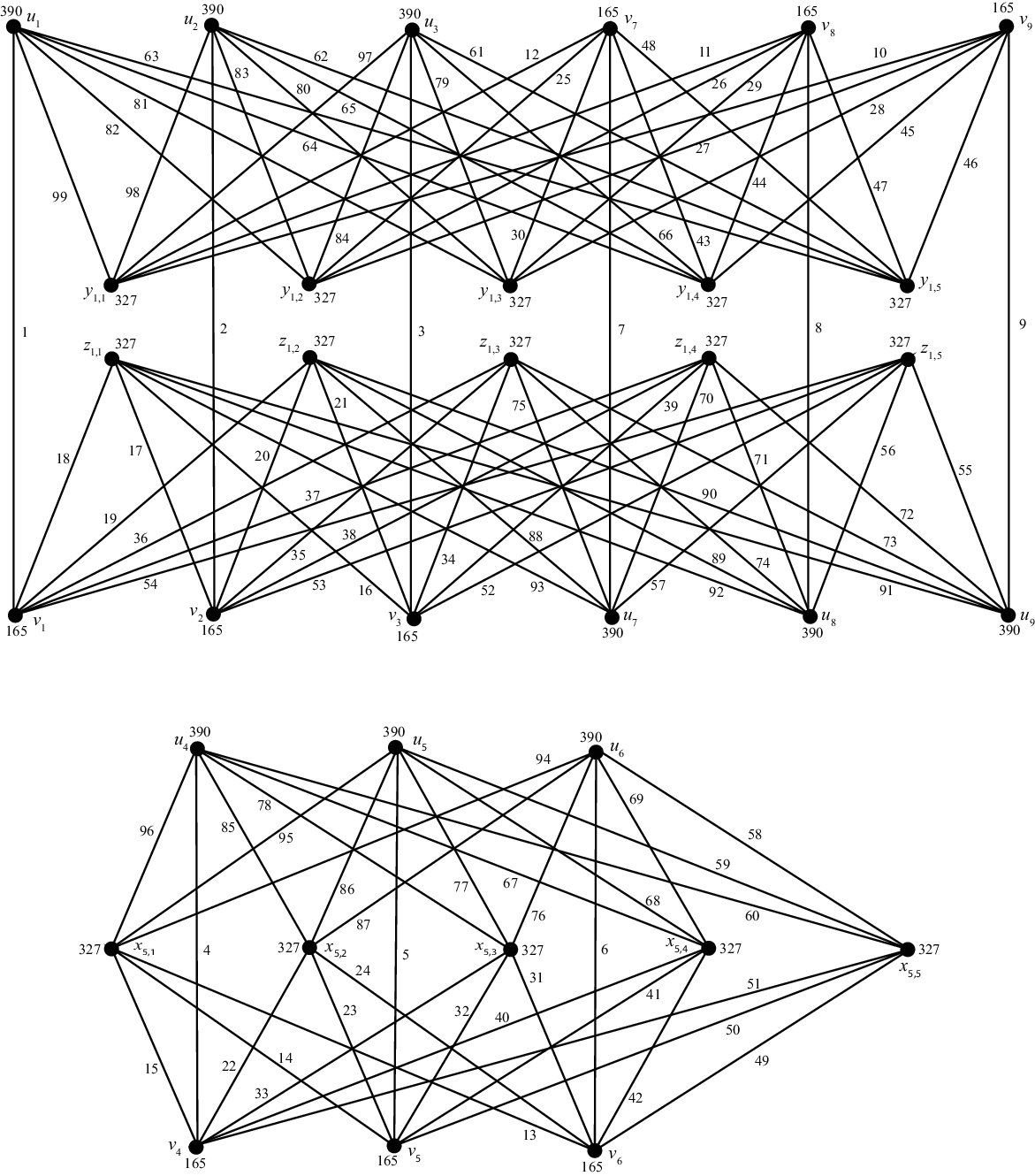, width= 15cm}}
\caption{$G_5(3,3)$ is a 6-regular tripartite graph.}\label{fig:G5(33))}
\end{figure}
\end{example}

\nt Note that we may also apply the delete-add process that gives us Theorem~2.6 in~\cite{LauShiu-join} to the graphs $G_{2n}(2r+1,2s+1)$ and $G_{2n+1}(2r+1,2s+1)$ to obtain two new families of (possibly connected or regular) tripartite graphs with local antimagic chromatic number 3. Denote the respective families of graph as $\mathcal R_{2n}(2r+1,2s+1)$ and $\mathcal R_{2n+1}(2r+1,2s+1)$. For example, from graph $G_4(3,3)$, we may remove the edges $v_9y_{1,1}$, $u_1y_{1,1}$ with labels $13,78$  and $u_4x_{5,1}, u_6x_{5,1}$ with labels $81,10$ respectively; and add the edges $v_9x_{5,1}$ with label 13, $u_1x_{5,1}$ with label 78, $u_4y_{1,1}$ with label 81, and $u_6y_{1,1}$ with label 10. The new graph is in $\mathcal R_4(3,3)$ and is connected. If we apply this process to $G_5(3,3)$ involving the edges with labels $99,10$ and $96,13$ respectively, we get a connected 6-regular graph in $\mathcal R_5(3,3)$. Thus, we have the following corollary with the proof omitted.

\begin{corollary} For $n, r, s\ge 1$, if $n=2s$, $\mathcal R_{2n+1}(2r+1,2s+1)$ is a family of (possibly connected) $(2n+2)$-regular tripartite graphs with local antimagic chromatic number $3$.
\end{corollary}

\section{Conclusions and Discussion}

In this paper, we constructed severy families of infinitely many tripartite graphs of size $(4n+1)\times (2k+1)$ and $(4n+3)\times (2k+1)$ respectively. We then use matrices to show that these graphs have local antimagic chromatic number 3. As a natural extension, we shall in another paper show that such families of graphs of size $(4n+1)\times 2k$ and $(4n+3)\times 2k$ respectively are bipartite but they also have local antimagic chromatic number 3. Interested readers may refer to~\cite{LSNP-even} for more related results.

\end{document}